\documentclass{article}

\newtheorem{theorem}{Theorem}[section]
\newtheorem{corollary}[theorem]{Corollary}

\newtheorem{lemma}[theorem]{Lemma}
\newtheorem{proposition}[theorem]{Proposition}

\evensidemargin\oddsidemargin
\begin{document}

\author{Vadim E. Levit and Eugen Mandrescu \\
Department of Computer Science\\
Holon Academic Institute of Technology\\
52 Golomb Str., P.O. Box 305\\
Holon 58102, ISRAEL}
\title{On $\alpha $-Critical Edges in K\"{o}nig-Egerv\'{a}ry Graphs}
\date{}
\maketitle

\begin{abstract}
The \textit{stability number} of a graph $G$, denoted by $\alpha (G)$, is
the cardinality of a stable set of maximum size in $G$. If $\alpha
(G-e)>\alpha (G)$, then $e$ is an $\alpha $\textit{-critical edge}, and if $%
\mu (G-e)<\mu (G)$, then $e$ is a $\mu $\textit{-critical edge}, where $\mu
(G)$ is the cardinality of a maximum matching in $G$. $G$ is a \textit{%
K\"{o}nig-Egerv\'{a}ry graph} if its order equals $\alpha (G)+\mu (G)$.
Beineke, Harary and Plummer have shown that the set of $\alpha $-critical
edges of a bipartite graph is a matching. In this paper we generalize this
statement to K\"{o}nig-Egerv\'{a}ry graphs. We also prove that in a
K\"{o}nig-Egerv\'{a}ry graph $\alpha $-critical edges are also $\mu $%
-critical, and that they coincide in bipartite graphs. Eventually, we deduce
that $\alpha (T)=\xi (T)+\eta (T)$ holds for any tree $T$, and characterize
the K\"{o}nig-Egerv\'{a}ry graphs enjoying this property, where $\xi (G)$ is
the number of $\alpha $-critical vertices of $G,$ and $\eta (G)$ is the
number of $\alpha $-critical edges of $G$.
\end{abstract}

\section{Introduction}

Throughout this paper $G=(V,E)$ is a simple (i.e., a finite, undirected,
loopless and without multiple edges) graph with vertex set $V=V(G)$, edge
set $E=E(G)$, and order $n(G)=|V(G)|$. If $X\subset V$, then $G[X]$ is the
subgraph of $G$ spanned by $X$. By $G-W$ we mean the subgraph $G[V-W]$ , if $%
W\subset V(G)$. For $F\subset E(G)$, by $G-F$ we denote the partial subgraph
of $G$ obtained by deleting the edges of $F$, and we use $G-e$, if $W$ $%
=\{e\}$. If $A,B$ $\subset V$ and $A\cap B=\emptyset $, then $(A,B)$ stands
for the set $\{e=ab:a\in A,b\in B,e\in E\}$. The neighborhood of a vertex $%
v\in V$ is the set $N(v)=\{w:w\in V$ \ \textit{and} $vw\in E\}$, and $%
N(A)=\cup \{N(v):v\in A\}$, $N[A]=A\cup N(A)$ for $A\subset V$.

A set $S$ of vertices is \textit{stable} if no two vertices from $S$ are
adjacent. A stable set of maximum size will be referred to as a \textit{%
maximum stable set} of $G$. The \textit{stability number }of $G$, denoted by 
$\alpha (G)$, is the cardinality of a maximum stable set\textit{\ }of $G$.
Let $\Omega (G)$ denotes the set $\{S:S$ \textit{is a maximum stable set of} 
$G\}$, $\sigma (G)=\left| \cap \{V-S:S\in \Omega (G)\}\right| $ and $\xi
(G)=\left| core(G)\right| $, where $core(G)=\cap \{S:S\in \Omega (G)\}$, 
\cite{levm3}. In other words, $\xi (G)$ equals the number of $\alpha $%
-critical vertices of $G$, (a vertex $v\in V(G)$ is $\alpha $\textit{%
-critical} provided $\alpha (G-v)<\alpha (G)$).

By $P_{n},C_{n},K_{n}$ we mean the chordless path on $n\geq 3$, the
chordless cycle on $n\geq $ $4$ vertices, and respectively the complete
graph on $n\geq 1$ vertices.

A matching (i.e., a set of non-incident edges of $G$) of maximum cardinality 
$\mu (G)$ is a \textit{maximum matching}, and a \textit{perfect matching} is
one covering all vertices of $G$. An edge $e\in E(G)$ is $\mu $-\textit{%
critical }provided $\mu (G-e)<\mu (G)$. By their definition, $\mu $-critical
edges of $G$ belong to all maximum matchings of $G$.

If $\alpha (G)+\mu (G)=n(G)$, then $G$ is called a \textit{%
K\"{o}nig-Egerv\'{a}ry graph}, \cite{dem}, \cite{ster}.\textit{\ }Properties
of these graphs were presented in several papers, like of Sterboul \cite
{ster}, Deming \cite{dem}, Lov\'{a}sz and Plummer \cite{lovpl}, Korach \cite
{kora}, Bourjolly and Pulleyblank \cite{bourpull}, Paschos and Demange \cite
{pasdema}, Levit and Mandrescu \cite{levm2}, \cite{levm4}. It is worth
observing that a disconnected graph is of K\"{o}nig-Egerv\'{a}ry type if and
only if all its connected components are K\"{o}nig-Egerv\'{a}ry graphs. In
this paper, by ''graph'' we mean a connected graph having at least one edge.

An edge $e\in E(G)$ is $\alpha $-\textit{critical} whenever $\alpha
(G-e)>\alpha (G)$. Let denote by $\eta (G)$ the number of $\alpha $-critical
edges of $G$. Notice that there are graphs in which: ($\mathit{a}$) any edge
is $\alpha $-critical (so-called $\alpha $\textit{-critical graphs}); e.g.,
all $C_{2n+1}$ for $n\geq 3$; ($\mathit{b}$) no edge is $\alpha $-critical;
e.g., all $C_{2n}$ for $n\geq 2$. More generally, Haynes et al., \cite{hayn}%
, have proved that a graph $G$ has no $\alpha $-critical edge if and only if 
$\left| N(x)\cap S\right| \geq 2$ holds for any $S\in \Omega (G)$ and every $%
x\in V(G)-S$.

Beineke, Harary and Plummer, \cite{BeiHarPlum}, have shown that any two
incident $\alpha $-critical edges of a graph lie on an odd cycle, and hence,
they deduce that no two $\alpha $-critical edges of a bipartite graph can
have a common endpoint. Independently, Zito, \cite{Zito}, has proved the
same result for trees using a different technique. Some variations and
strengthenings of these results are discussed in \cite{Suranyi}, \cite
{Wessel}, and \cite{Toft}.

In this paper we generalize the above assertion to K\"{o}nig-Egerv\'{a}ry
graphs. We also show that $\alpha $-critical edges are $\mu $-critical in a
K\"{o}nig-Egerv\'{a}ry graph, and that they coincide in bipartite graphs. As
a corollary, we obtain one result of Zito, \cite{Zito}, stating that a
vertex $v$ is in some but not in all maximum stable sets of a tree $T$ if
and only if $v$ is an endpoint of an $\alpha $-critical edge of $T$. In the
sequel, we analyze other relationships between $\alpha $-critical edges and $%
\mu $-critical edges in a K\"{o}nig-Egerv\'{a}ry graph, and its
corresponding implications to equalities and inequalities linking $\alpha
(G) $, $\xi (G)$, $\eta (G)$, $\sigma (G)$, and $\mu (G)$. Eventually, we
infer that $\alpha (T)=\xi (T)+\eta (T),\sigma (T)+\eta (T)=\mu (T)$ and $%
\xi (T)+2\eta (T)+\sigma (T)=n(T)$ holds for any tree $T$, and characterize
the K\"{o}nig-Egerv\'{a}ry graphs having these properties.

\section{$\alpha $-Critical and $\mu $-Critical Edges}

\begin{figure}[h]
\setlength{\unitlength}{1cm}%
\begin{picture}(5,1)\thicklines

  \multiput(6,0)(1,0){3}{\circle*{0.29}}
  \put(7,1){\circle*{0.29}}
  \put(6,0){\line(1,0){2}}
  \put(7,0){\line(0,1){1}} 
  \put(8,0){\line(-1,1){1}} 
 \put(6.5,0.25){\makebox(0,0){$e$}} 
 \end{picture}
\caption{Graph\emph{\ }$K_{3}+e.$}
\label{fig1}
\end{figure}
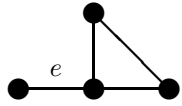

According to a well-known result of K\"{o}nig, \cite{koen}, and
Egerv\'{a}ry, \cite{eger}, any bipartite graph is a K\"{o}nig-Egerv\'{a}ry
graph. It is easy to see that this class includes also some non-bipartite
graphs (see, for instance, the graph $K_{3}+e$ in Figure \ref{fig1}).

If $G_{i}=(V_{i},E_{i}),i=1,2$, are two disjoint graphs, then $G=G_{1}*G_{2}$
is defined as the graph with $V(G)=V(G_{1})\cup V(G_{2})$, and 
\[
E(G)=E(G_{1})\cup E(G_{2})\cup \{xy:for\ some\ x\in V(G_{1})\ and\ y\in
V(G_{2})\}. 
\]

Clearly, if $H_{1},H_{2}$ are subgraphs of a graph $G$ such that $%
V(G)=V(H_{1})\cup V(H_{2})$ and $V(H_{1})\cap V(H_{2})=$ $\emptyset $, then $%
G=H_{1}*H_{2}$, i.e., any graph of order at least two admits such
decompositions. However, some particular cases are of special interest. For
instance, if: $E(H_{i})=\emptyset ,i=1,2,$ then $G=H_{1}*H_{2}$ is
bipartite; $E(H_{1})=\emptyset $ and $H_{2}$ is complete, then $%
G=H_{1}*H_{2} $ is a \textit{split graph }\cite{hmmr}.

The following result shows that the K\"{o}nig-Egerv\'{a}ry graphs are, in
this sense, between these two ''extreme'' situations. The equivalence of the
first and the third parts of this proposition was proposed by Klee and
included in \cite{lawl} without proof (private communication).

\begin{proposition}
\label{prop1}\cite{levm4} The following assertions are equivalent:

($\mathit{i}$) $G$ is a K\"{o}nig-Egerv\'{a}ry graph;

($\mathit{ii}$) $G=H_{1}*H_{2}$, where $V(H_{1})=S\in \Omega (G)$ and $%
n(H_{1})\geq \mu (G)=n(H_{2})$;

($\mathit{iii}$) $G=H_{1}*H_{2}$, where $V(H_{1})=S$ is a stable set in $G,$ 
$\left| S\right| \geq n(H_{2})$ and $(S,V(H_{2}))$ contains a matching $M$
with $\left| M\right| =n(H_{2})$.
\end{proposition}

In the sequel, we shall often represent a K\"{o}nig-Egerv\'{a}ry graph $G$
as $G=S*H$, where $S\in \Omega (G)$ and $H=G[V-S]$ has $n(H)=\mu (G)$.

\begin{lemma}
\label{lem2}\cite{levm4} If $G=(V,E)$ is a K\"{o}nig-Egerv\'{a}ry graph,
then any maximum matching of $G$ is contained in $(S,V-S)$, where $S\in
\Omega (G)$.
\end{lemma}

Clearly, Lemma \ref{lem2} is not valid for any graph. For instance, $K_{4}$
is a counterexample. Moreover, $K_{4}$ has $\alpha $-critical edges that are
incident. Nevertheless, there are graphs having only non-incident $\alpha $%
-critical edges.

\begin{theorem}
\label{TH1}If $G$ is a K\"{o}nig-Egerv\'{a}ry graph, then the following
assertions hold:

($\mathit{i}$) for any $\alpha $-critical edge $e$ of $G$, the graph $G-e$
is still a K\"{o}nig-Egerv\'{a}ry graph;

($\mathit{ii}$) any $\alpha $-critical edge of $G$ is also $\mu $-critical;

($\mathit{iii}$) the $\alpha $-critical edges of $G$ form a matching.
\end{theorem}

\setlength {\parindent}{0.0cm}\textbf{Proof.} ($\mathit{i}$) If $e=xy$ is an 
$\alpha $-critical edge $G$, then there is some $S\in \Omega (G)$ such that
either $N(x)\cap S=\{y\}$ or $N(y)\cap S=\{x\}$. Suppose that $y\in S$.
Since $S\in \Omega (G)$, we get, by Proposition \ref{prop1}, that $G=S*H$,
where $H=G[V-S]$ has $\mu (G)=n(H)=\left| M\right| $ and $M$ is a maximum
matching of $G$, included, by Lemma \ref{lem2}, in $(S,V(G)-S)$. Hence, it
follows that $G-e=S^{\prime }*V(H^{\prime })$, where $S^{\prime }=S\cup
\{x\}\in \Omega (G-e)$ and $n(H^{\prime })=\left| M-\{e\}\right| $.
According to Proposition \ref{prop1}($\mathit{iii}$), we infer that $G-e$ is
also a K\"{o}nig-Egerv\'{a}ry graph.\setlength
{\parindent}{3.45ex}

($\mathit{ii}$) If $e\in E(G)$ is an $\alpha $-critical edge of $G$, then
according to ($\mathit{i}$) we obtain: 
\[
n(G)=\alpha (G)+\mu (G)\leq \alpha (G-e)+\mu (G-e)=\alpha (G)+1+\mu
(G-e)=n(G-e), 
\]
and this implies $\mu (G)=1+\mu (G-e)$, i.e., $e$ is also $\mu $-critical.

($\mathit{iii}$) Let $e_{1},e_{2}$ be two $\alpha $-critical edges of $G$.
We have to show that they are not incident. According to second part ($%
\mathit{ii}$), both edges are also $\mu $-critical. Hence, it follows that $%
e_{1},e_{2}\in \cap \{M:M\ is\ a\ maximum\ matching\ of\ G\}$ and this
ensures that $e_{1},e_{2}$ have no common endpoint. Consequently, the set of
all $\alpha $-critical edges of $G$ yields a matching. \rule{2mm}{2mm}%
\newline

Notice that:

($\mathit{a}$) Theorem \ref{TH1}($\mathit{i}$) is not true for any $\mu $%
-critical edge of a K\"{o}nig-Egerv\'{a}ry graph; e.g., the edge $e$ of $%
G=K_{3}+e$ is $\mu $-critical, but $G-e$ is not a K\"{o}nig-Egerv\'{a}ry
graph;

($\mathit{b}$) Theorem \ref{TH1}($\mathit{ii}$) is not true for any graph;
e.g., all the edges of $K_{3}$ are $\alpha $-critical, but none is also $\mu 
$-critical;

($\mathit{c}$) the converse of Theorem \ref{TH1}($\mathit{ii}$) is not valid
for any K\"{o}nig-Egerv\'{a}ry graph; e.g., the edge $e$ of graph $K_{3}+e$
is $\mu $-critical, but is not also $\alpha $-critical. However, as we shall
see later, (namely Proposition \ref{prop3}), the $\mu $-critical edges are
also $\alpha $-critical in the case of bipartite graphs.

\begin{corollary}
A K\"{o}nig-Egerv\'{a}ry graph is $\alpha $-critical if and only if it is
isomorphic to $K_{2}$.
\end{corollary}

Since any bipartite graph is also a K\"{o}nig-Egerv\'{a}ry graph, we obtain
the following statement, due to Beineke, Harary and Plummer.

\begin{theorem}
\cite{BeiHarPlum} No two $\alpha $-critical edges of a bipartite graph are
incident.
\end{theorem}

\begin{proposition}
\label{prop3}If $G$ is a bipartite graph, then its $\alpha $-critical edges
coincide with its $\mu $-critical edges.
\end{proposition}

\setlength {\parindent}{0.0cm}\textbf{Proof.} By Theorem \ref{TH1}($\mathit{%
ii}$), it suffices to show that any $\mu $-critical edge $e$ of $G$ is also $%
\alpha $-critical. Since $G-e$ is still bipartite, and hence, also a
K\"{o}nig-Egerv\'{a}ry graph, it follows that $\alpha (G-e)+\mu
(G-e)=n(G)=\alpha (G)+\mu (G)=\alpha (G)+1+\mu (G-e)$, and this implies $%
\alpha (G-e)>\alpha (G)$, i.e., $e$ is an $\alpha $-critical edge of $G$. 
\rule{2mm}{2mm}\setlength
{\parindent}{3.45ex}\newline

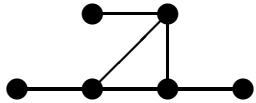
\begin{figure}[h]
\setlength{\unitlength}{1cm}%
\begin{picture}(5,1)\thicklines

  \multiput(5,0)(1,0){4}{\circle*{0.29}}
  \put(7,1){\circle*{0.29}}
  \put(6,1){\circle*{0.29}}
  \put(5,0){\line(1,0){3}}
  \put(6,1){\line(1,0){1}}
  \put(7,0){\line(0,1){1}} 
  \put(6,0){\line(1,1){1}} 
  
 \end{picture}
\caption{A Koenig-Egervary graph whose all $\mu $-critical edges are $\alpha 
$-critical.}
\label{fig4}
\end{figure}

In Theorem \ref{Th2} we will meet another type of K\"{o}nig-Egerv\'{a}ry
graphs with this property. Notice that there are also non-bipartite
K\"{o}nig-Egerv\'{a}ry graphs in which their $\mu $-critical edges are $%
\alpha $-critical (see the graph in Figure \ref{fig4}).

It is well-known that if a tree has a perfect matching, then it is unique.
Consequently, we obtain:

\begin{corollary}
\label{cor4}A tree has a perfect matching if and only if the set of its $%
\alpha $-critical edges forms a maximal matching of the tree.
\end{corollary}

Using the definition of K\"{o}nig-Egerv\'{a}ry graphs and the fact that $\mu
(G)\leq n(G)/2$ is true for any graph $G$, we get:

\begin{lemma}
\label{lem1}If $G$ admits a perfect matching, then $G$ is a
K\"{o}nig-Egerv\'{a}ry graph if and only if $\alpha (G)=\mu (G)$. If $G$ is
a K\"{o}nig-Egerv\'{a}ry graph, then $\mu (G)\leq \alpha (G)$.
\end{lemma}

Combining Corollary \ref{cor4} and Lemma \ref{lem1}, we get the following
result from \cite{Zito}.

\begin{corollary}
\cite{Zito}\label{cor3} If a tree $T$ has a perfect matching $M$, then all
the edges of $M$ are $\alpha $-critical and $2\alpha (T)=n(T)$.
\end{corollary}

\begin{proposition}
\label{prop5}If $G=(V,E)$ is a K\"{o}nig-Egerv\'{a}ry graph, then the
following assertions are true:

($\mathit{i}$) any $S\in \Omega (G)$ meets each $\mu $-critical edge in
exactly one vertex;

($\mathit{ii}$) any $S\in \Omega (G)$ meets each $\alpha $-critical edge in
exactly one vertex;

($\mathit{iii}$) if $G$ has a maximal matching consisting of only $\alpha $%
-critical edges, then it is the unique perfect matching of $G$.
\end{proposition}

\setlength {\parindent}{0.0cm}\textbf{Proof.} ($\mathit{i}$) $\mathit{and}$ (%
$\mathit{ii}$) By Theorem \ref{TH1}($\mathit{ii}$), any $\alpha $-critical
edge of $G$ is also $\mu $-critical. Consequently, we infer that 
\[
\{e\in E:e\ is\ \alpha -critical\}\subseteq \cap \{M:M\ is\ a\ maximum\
matching\ of\ G\}\subseteq (S,V-S) 
\]
holds for any $S\in \Omega (G)$, according to Lemma \ref{lem2}. It follows
that if $e=xy$ is an $\alpha $-critical or a $\mu $-critical edge of $G$,
then any $S\in \Omega (G)$ contains one of $x$ and $y$, (since clearly, no
stable set may contain both $x$ and $y$).\setlength
{\parindent}{3.45ex}

($\mathit{iii}$) Let $M$ be a maximal matching of $G$ consisting of only $%
\alpha $-critical edges. By Theorem \ref{TH1}, all the edges of $M$ are also 
$\mu $-critical. Therefore, we infer that $M$ is included in any maximum
matching of $G$, and because $M$ is a maximal matching, it results that $M$
is the unique maximum matching of $G$. Suppose, on the contrary, that $M$ is
not perfect, and let $S\in \Omega (G)$. According to Proposition \ref{prop1}%
, $G$ can be written as $G=S*H$, with $n(H)=\left| M\right| =\mu (G)$, and
by Lemma \ref{lem2} we have that $M\subseteq (S,V-S)$. Since $G$ is a
K\"{o}nig-Egerv\'{a}ry graph without perfect matchings, Lemma \ref{lem1}
implies $\left| S\right| =\alpha (G)>\mu (G)=\left| M\right| $. Hence, it
follows that there are at least two vertices $v_{1},v_{2}\in S$ having a
common neighbor $w\in V(H)$ and such that one of them, say $v_{1}$, is
unmatched by $M$ and $v_{2}w\in M$. Thus, $M\cup \{v_{1}w\}-\{v_{2}w\}$ is
another maximum matching of $G$, in contradiction with the uniqueness of $M$%
. Consequently, $M$ must be also perfect. \rule{2mm}{2mm}\newline

For trees, Proposition \ref{prop5}($\mathit{ii}$) was proved by Zito in \cite
{Zito}.

Notice that the matching in Proposition \ref{prop5}($\mathit{iii}$) is not
necessarily formed by pendant edges; e.g., $P_{6}$ has such a matching.
Concerning the uniqueness of this matching, it is worth mentioning that: ($%
\mathit{a}$) if $G$ is not a K\"{o}nig-Egerv\'{a}ry graph, then it may have
several different maximum matchings consisting of only $\alpha $-critical
edges (e.g., $C_{5}$); ($\mathit{b}$) if a K\"{o}nig-Egerv\'{a}ry graph has
a unique perfect matching, then it may contain non-$\alpha $-critical edges
(e.g., the edge $e$ of $K_{3}+e$ is not $\alpha $-critical, but it belongs
to the unique perfect matching of $K_{3}+e$).

\section{Equalities and Inequalities between Parameters}

If $v\in N(core(G))$, then clearly follows that $v\in V(G)-S$, for any $S\in
\Omega (G)$, that is $N(core(G))\subseteq \cap \{V-S:S\in \Omega (G)\}$
holds for any graph $G$.

\begin{lemma}
\label{prop11}\cite{levm4} If $G=(V,E)$ is a K\"{o}nig-Egerv\'{a}ry graph,
then

$N(core(G))=\cap \{V-S:S\in \Omega (G)\}.$
\end{lemma}

\begin{figure}[h]
\setlength{\unitlength}{1cm}%
\begin{picture}(5,1)\thicklines
  \multiput(3,0)(1,0){2}{\circle*{0.25}}
  \multiput(3,1)(1,0){2}{\circle*{0.25}}
  \put(3.5,0.5){\circle*{0.25}}
  \put(3,0){\line(1,0){1}}
  \put(3,1){\line(1,0){1}}
  \put(3,0){\line(0,1){1}}
  \put(4,0){\line(0,1){1}}
  \put(3,0){\line(1,1){1}} 
  \put(3,1){\line(1,-1){1}} 
  \put(3.5,0.8){\makebox(0,0){$v$}}
  \put(2,0.5){\makebox(0,0){$(a)$}}

  \multiput(7,0)(1,0){2}{\circle*{0.25}}
  \multiput(7,1)(1,0){2}{\circle*{0.25}}
  \put(8.5,0.5){\circle*{0.25}}
  \put(7,0){\line(1,0){1}}
  \put(7,1){\line(1,0){1}}
  \put(7,0){\line(3,1){1.5}}
  \put(7,1){\line(3,-1){1.5}}
  \put(8,0){\line(0,1){1}}
  \put(7,0){\line(1,1){1}} 
  \put(7,1){\line(1,-1){1}} 
  \put(8,0){\line(1,1){0.5}}
  \put(8,1){\line(1,-1){0.5}} 
  \put(6,0.5){\makebox(0,0){$(b)$}}
\end{picture}
\caption{(a) $G$ is non-K\"{o}nig-Egerv\'{a}ry with $N(core(G))\neq \cap
\{V-S:S\in \Omega (G)\}$; (b) $G$ is a non-K\"{o}nig-Egerv\'{a}ry graph with 
$N(core(G))=\cap \{V-S:S\in \Omega (G)\}$.}
\label{fig100}
\end{figure}
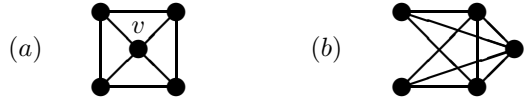

Notice that there are graphs that do not enjoy the above equality, for
example, the graph $G$ in Figure \ref{fig100}($\mathit{a}$) has $%
N(core(G))=\emptyset $ and $\cap \{V-S:S\in \Omega (G)\}=\{v\}$. There exist
non-K\"{o}nig-Egerv\'{a}ry graphs for which $N(core(G))=\cap \{V-S:S\in
\Omega (G)\}$, (see, for instance, the graph $G$ from Figure \ref{fig100}($%
\mathit{b}$)).

\begin{proposition}
\label{prop9}If $G=(V,E)$ is a K\"{o}nig-Egerv\'{a}ry graph, $%
G_{0}=G-N[core(G)]$ and $S\in \Omega (G)$, then the following assertions are
true:

($\mathit{i}$) $\left| core(G)\right| \geq \left| N(core(G))\right| $;

($\mathit{ii}$) $\left| S-core(G)\right| =\left| V-S-N(core(G))\right| $;

($\mathit{iii}$) $G_{0}$ has a perfect matching and it is also a
K\"{o}nig-Egerv\'{a}ry graph.
\end{proposition}

\setlength {\parindent}{0.0cm}\textbf{Proof.} According to Proposition \ref
{prop1}, $G$ can be written as $G=S*H$, where $H=G[V-S]$ has $n(H)=\mu (G)$.
Let denote $A=S-core(G)$ and $B=V(H)-N(core(G))$. In \cite{levm3} it has
been proved that $\left| A\right| \leq \left| B\right| $ holds for any graph 
$G$. Since $\cap \{V-S:S\in \Omega (G)\}\subseteq V(H)$, and $%
N(core(G))=\cap \{V-S:S\in \Omega (G)\}$ (see Lemma \ref{prop11}), we obtain 
$B=V(H)-\cap \{V-S:S\in \Omega (G)\}$.\setlength
{\parindent}{3.45ex}

($\mathit{i}$) Since $\left| A\right| +\left| core(G)\right| =\alpha (G)\geq
\mu (G)=n(H)=\left| B\right| +\left| N(core(G)\right| $ and, on the other
hand $\left| A\right| \leq \left| B\right| $, it follows that $\left|
core(G)\right| \geq \left| N(core(G))\right| $.

($\mathit{ii}$) Let $M$ be a maximum matching in $G$. Since $G$ is a
K\"{o}nig-Egerv\'{a}ry graph, Lemma \ref{lem2} ensures that $M$ is included
in $(S,V(H))$, and $\left| M\right| =\mu (G)=n(H)$. The matching $M$ matches 
$B$ into $A$, because there are no edges connecting $B$ and $core(G)$.
Hence, $\left| B\right| \leq \left| A\right| $. Together with $\left|
A\right| \leq \left| B\right| $ $\left| A\right| \leq \left| B\right| $, it
implies $\left| A\right| =\left| B\right| $, i.e., $\left| S-core(G)\right|
=\left| V-S-N(core(G))\right| $, and that $M\cap (A,B)$ is a perfect
matching of $G[A\cup B]$.

($\mathit{iii}$) Since, in fact, $G_{0}=G[A\cup B]$, it follows necessarily
that $G_{0}$ has a perfect matching. In addition, because $A$ is stable, we
get $\alpha (G_{0})\leq \mu (G_{0})=\left| A\right| \leq \alpha (G_{0})$,
i.e., $\alpha (G_{0})=\mu (G_{0})$, and according to Lemma \ref{lem1}, $%
G_{0} $ must be also a K\"{o}nig-Egerv\'{a}ry graph. \rule{2mm}{2mm}

\begin{corollary}
\label{cor2}If $G$ is a K\"{o}nig-Egerv\'{a}ry graph, then $\alpha
(G)+\sigma (G)=\mu (G)+\xi (G)$.
\end{corollary}

\setlength {\parindent}{0.0cm}\textbf{Proof.} By Lemma \ref{prop11}, $%
N(core(G))=\cap \{V-S:S\in \Omega (G)\}$ and according to Proposition \ref
{prop9}($\mathit{ii}$), $\left| S-core(G)\right| =\left|
V-S-N(core(G))\right| $. Hence, we obtain that $\alpha (G)-\xi (G)=\left|
S-core(G)\right| =\left| V-S-N(core(G))\right| =\mu (G)-\sigma (G)$. \rule%
{2mm}{2mm}\setlength
{\parindent}{3.45ex}\newline

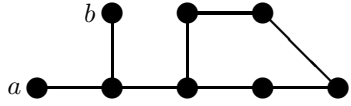
\begin{figure}[h]
\setlength{\unitlength}{1cm}%
\begin{picture}(5,1.5)\thicklines

  \multiput(4,0)(1,0){5}{\circle*{0.29}}
  \multiput(5,1)(1,0){3}{\circle*{0.29}}
  \put(4,0){\line(1,0){4}}
  \put(5,0){\line(0,1){1}} 
  \put(6,0){\line(0,1){1}} 
  \put(6,1){\line(1,0){1}}
  \put(7,1){\line(1,-1){1}}

\put(3.7,0){\makebox(0,0){$a$}}
\put(4.7,1){\makebox(0,0){$b$}}

 \end{picture}
\caption{$G$ is a non-K\"{o}nig-Egervary graph with $core(G)=\{a,b\}$ and $%
G_{0}=C_{5}$.}
\label{fig7}
\end{figure}

Let us observe that there exist non-K\"{o}nig-Egerv\'{a}ry graphs satisfying
the equality $\alpha (G)+\sigma (G)=\mu (G)+\xi (G)$ (see graph $W_{1}$ in
Figure \ref{fig2}). It is also interesting to notice that there exists a
non-K\"{o}nig-Egerv\'{a}ry graph enjoying the property that its subgraph $%
G_{0}=G-N[core(G)]$ has a perfect matching (see Figure \ref{fig5}). Figure 
\ref{fig7} shows a non-K\"{o}nig-Egerv\'{a}ry graph $G$ whose $G_{0}$ has no
perfect matching.

\begin{lemma}
\label{prop6}Let $G=(V,E)$ and $G_{0}=G-N[core(G)]$. Then the following
assertions are valid:

($\mathit{i}$) no $\alpha $-critical edge in $G$ has an endpoint in $%
N[core(G)]$;

($\mathit{ii}$) $\alpha (G)=\alpha (G_{0})+\xi (G),\Omega (G_{0})=\{S\cap
V(G_{0}):S\in \Omega (G)\},core(G_{0})=\emptyset $;

($\mathit{iii}$) $e=xy$ is an $\alpha $-critical edge of $G$ if and only if $%
e$ is an $\alpha $-critical edge of $G_{0}$.
\end{lemma}

\setlength {\parindent}{0.0cm}\textbf{Proof.} ($\mathit{i}$) Let $e=xy$ be
an $\alpha $-critical edge in $G$, and let $\overline{S}\in \Omega (G-e)$.
Since $\left| \overline{S}\right| =\alpha (G-e)>\alpha (G)$, it follows that 
$x,y\in \overline{S}$ and $\overline{S}-\{x\},\overline{S}-\{y\}\in \Omega
(G)$. Now, the inclusion $N(core(G))\subseteq \cap \{V-S:S\in \Omega (G)\}$
completes the proof that no $\alpha $-critical edge in $G$ has an endpoint
in $N(core(G))$, and respectively, in $core(G)$.%
\setlength
{\parindent}{3.45ex}

($\mathit{ii}$) By definition of $G_{0}$, if $S\in \Omega (G)$, then $%
S-core(G)=S\cap V(G_{0})$, and therefore 
\[
\alpha (G)-\xi (G)=\left| S-core(G)\right| \leq \alpha (G_{0}). 
\]
For any $S_{G_{0}}\in \Omega (G_{0})$ we have that $S_{G_{0}}\cup core(G)$
is stable, and hence 
\[
\left| S_{G_{0}}\cup core(G)\right| =\alpha (G_{0})+\xi (G)\leq \alpha (G). 
\]
Consequently, we get $\alpha (G)=\alpha (G_{0})+\xi (G)$. Now it is easy to
check that $\Omega (G_{0})=\{S\cap V(G_{0}):S\in \Omega (G)\}$ and $%
core(G_{0})=\emptyset $.

($\mathit{iii}$) Let $e=xy$ be an $\alpha $-critical edge of $G$. By ($%
\mathit{i}$), we infer that $e\in E(G_{0})$, and as we saw above, there is
some stable set $S_{xy}$ such that $S_{xy}\cup \{x\},S_{xy}\cup \{y\}\in
\Omega (G)$ and $S_{xy}\cup \{x,y\}\in \Omega (G-e)$. Hence, ($\mathit{ii}$)
implies that 
\[
V(G_{0})\cap (S_{xy}\cup \{x\}),V(G_{0})\cap (S_{xy}\cup \{y\})\in \Omega
(G_{0})\ and\ V(G_{0})\cap (S_{xy}\cup \{x,y\})\in \Omega (G_{0}-e), 
\]
because $V(G_{0})\cap (S_{xy}\cup \{x,y\})$ is stable in $G_{0}-e$ and
larger than $V(G_{0})\cap (S_{xy}\cup \{x\})$. Therefore, $e$ is $\alpha $%
-critical in $G_{0}$, as well. Similarly, we can show that any $\alpha $%
-critical edge of $G_{0}$ is $\alpha $-critical in $G$ too. \rule{2mm}{2mm}

\begin{proposition}
\label{prop7}If $G$ is a K\"{o}nig-Egerv\'{a}ry graph, then

($\mathit{i}$) $\xi (G)+\eta (G)\leq \alpha (G)$;

($\mathit{ii}$) $\sigma (G)+\eta (G)\leq \mu (G)$;

($\mathit{iii}$) $\xi (G)+2\eta (G)+\sigma (G)\leq n(G)$.
\end{proposition}

\setlength {\parindent}{0.0cm}\textbf{Proof.} For any $S\in \Omega (G)$, we
have that $core(G)\subseteq S$, and by Lemma \ref{prop6}($\mathit{i}$), no $%
\alpha $-critical edge has an endpoint in $core(G)$. In addition, according
to Proposition \ref{prop5}($\mathit{ii}$), $S$ meets each $\alpha $-critical
edge in exactly one vertex. Hence, it follows that $\xi (G)+\eta (G)\leq
\alpha (G)$, and using Corollary \ref{cor2} we obtain ($\mathit{ii}$).
Clearly, ($\mathit{iii}$) follows from ($\mathit{i}$) and ($\mathit{ii}$). 
\rule{2mm}{2mm}\setlength
{\parindent}{3.45ex}\newline

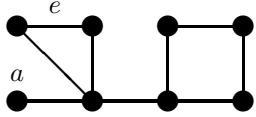
\begin{figure}[h]
\setlength{\unitlength}{1.0cm} 
\begin{picture}(5,1)\thicklines
  \multiput(5,0)(1,0){4}{\circle*{0.29}}
  \multiput(5,1)(1,0){4}{\circle*{0.29}}
  \put(5,0){\line(1,0){3}}
  \put(5,1){\line(1,0){1}}
  \put(5,1){\line(1,-1){1}}
  \put(7,1){\line(1,0){1}}
  \multiput(6,0)(1,0){3}{\line(0,1){1}}
  \put(5,0.35){\makebox(0,0){$a$}}
  \put(5.5,1.25){\makebox(0,0){$e$}}
 \end{picture}
\caption{${W}$ is a non-bipartite Koenig-Egervary graph and $\xi (W)+\eta
(W)<\alpha (W)$.}
\label{fig3}
\end{figure}

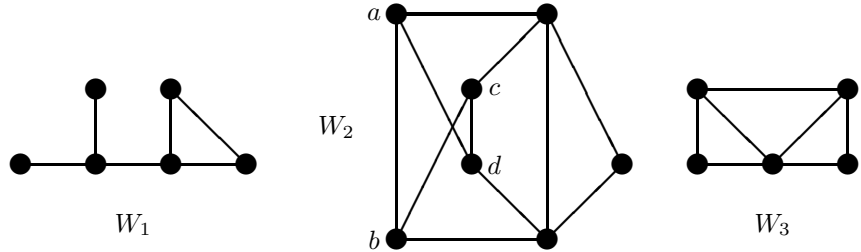
\begin{figure}[h]
\setlength{\unitlength}{1cm}%
\begin{picture}(5,3)\thicklines

  \multiput(1,1)(1,0){4}{\circle*{0.29}}
  \put(3,2){\circle*{0.29}}
  \put(2,2){\circle*{0.29}}
  \put(2,1){\line(0,1){1}}
  \put(1,1){\line(1,0){3}}
  \put(3,1){\line(0,1){1}} 
  \put(4,1){\line(-1,1){1}} 
\put(2.5,0.2){\makebox(0,0){$W_{1}$}} 
  
\multiput(6,0)(2,0){2}{\circle*{0.29}}
\multiput(6,3)(2,0){2}{\circle*{0.29}}
\put(7,1){\circle*{0.29}}
\put(7,2){\circle*{0.29}}
\put(9,1){\circle*{0.29}}
\put(6,0){\line(1,0){2}}
\put(6,3){\line(1,0){2}}
\put(6,0){\line(0,1){3}}
\put(8,0){\line(0,1){3}}
\put(6,0){\line(1,2){1}}
\put(6,3){\line(1,-2){1}}
\put(7,1){\line(0,1){1}}
\put(7,2){\line(1,1){1}}
\put(7,1){\line(1,-1){1}}
\put(8,0){\line(1,1){1}}
\put(8,3){\line(1,-2){1}}
\put(5.7,3){\makebox(0,0){$a$}} 
\put(5.7,0){\makebox(0,0){$b$}} 
\put(7.3,2){\makebox(0,0){$c$}} 
\put(7.3,1){\makebox(0,0){$d$}} 
\put(5.2,1.5){\makebox(0,0){$W_{2}$}} 

  \multiput(10,1)(1,0){3}{\circle*{0.29}}
  \multiput(10,2)(2,0){2}{\circle*{0.29}}
  \put(10,2){\line(1,0){2}}
  \put(10,1){\line(1,0){2}}
  \put(10,1){\line(0,1){1}}
  \put(12,1){\line(0,1){1}}
  \put(10,2){\line(1,-1){1}} 
  \put(11,1){\line(1,1){1}}
\put(11,0.2){\makebox(0,0){$W_{3}$}}
 \end{picture}
\caption{Non-Koenig-Egervary graphs.}
\label{fig2}
\end{figure}

Notice that $\xi (K_{3}+e)+\eta (K_{3}+e)=\alpha (K_{3}+e)$ and also $\eta
(K_{3}+e)+\sigma (K_{3}+e)=\mu (K_{3}+e)$, but there are
K\"{o}nig-Egerv\'{a}ry graphs satisfying $\xi (G)+\eta (G)<\alpha (G)$ and $%
\eta (G)+\sigma (G)<\mu (G)$. For instance, $G=C_{6}$, and also the graph $W$
in Figure \ref{fig3} is a K\"{o}nig-Egerv\'{a}ry non-bipartite graph that
has $\eta (W)=\left| \{e\}\right| =1,\xi (W)=\left| \{a\}\right| =1=\sigma
(W),\alpha (W)=\mu (W)=4$.

Observe that Proposition \ref{prop7} is not true for general graphs; e.g.,
the graph $W_{1}$ in Figure \ref{fig2} has $\alpha (W_{1})=3,\mu
(W_{1})=2,\eta (W_{1})=3,\xi (W_{1})=2,\sigma (W_{1})=1$. However, there are
non-K\"{o}nig-Egerv\'{a}ry graphs satisfying $\xi (G)+\eta (G)<\alpha (G)$
and $\eta (G)+\sigma (G)<\mu (G)$, for example, the graph $W_{2}$ in Figure 
\ref{fig2} has $\alpha (W_{2})=3,\eta (W_{2})=\left| \{ab,cd\}\right| ,\xi
(W_{2})=\sigma (W_{2})=0$. There also exist non-K\"{o}nig-Egerv\'{a}ry
graphs satisfying $\xi (G)+\eta (G)=\alpha (G)$ and $\eta (G)+\sigma (G)=\mu
(G)$, e.g., the graph $W_{3}$ in Figure \ref{fig2}. Nevertheless, $\xi
(K_{5}-e)+\eta (K_{5}-e)=\alpha (K_{5}-e)$, but $\eta (K_{5}-e)+\sigma
(K_{5}-e)>\mu (K_{5}-e)$.

\begin{proposition}
\label{prop10}If $G$ is a K\"{o}nig-Egerv\'{a}ry graph, then the following
assertions are equivalent:

($\mathit{i}$) $\xi (G)+\eta (G)=\alpha (G)$;

($\mathit{ii}$) $\sigma (G)+\eta (G)=\mu (G)$;

($\mathit{iii}$) $\xi (G)+2\eta (G)+\sigma (G)=n(G)$.
\end{proposition}

\setlength {\parindent}{0.0cm}\textbf{Proof.} Suppose that $\xi (G)+\eta
(G)=\alpha (G)$. According to Corollary \ref{cor2}, we get that $\mu
(G)=\alpha (G)+\sigma (G)-\xi (G)=\xi (G)+\eta (G)+\sigma (G)-\xi (G)=\eta
(G)+\sigma (G).$ The converse is proven in the same way.%
\setlength
{\parindent}{3.45ex}

Suppose $\xi (G)+2\eta (G)+\sigma (G)=n(G)$. Proposition \ref{prop7} claims
that $\xi (G)+\eta (G)\leq \alpha (G)$ and $\sigma (G)+\eta (G)\leq \mu (G)$%
. Together with $\alpha (G)+\mu (G)=n(G)$, which is true for
K\"{o}nig-Egerv\'{a}ry graphs, it gives us the two equalities needed.
Conversely, if, for instance, $\xi (G)+\eta (G)=\alpha (G)$ then, as we
already proved, $\sigma (G)+\eta (G)=\mu (G)$. Summing these two equalities
we obtain $\xi (G)+2\eta (G)+\sigma (G)=n(G)$. \rule{2mm}{2mm}

\section{K\"{o}nig-Egerv\'{a}ry Graphs for which $\xi +\eta =\alpha $}

\begin{lemma}
\label{lem3}Let $G$ be a K\"{o}nig-Egerv\'{a}ry graph and $G_{0}=G-N[core(G)]
$. If $G_{0}$ has a unique perfect matching then its $\alpha $-critical
edges coincide with its $\mu $-critical edges.
\end{lemma}

\setlength {\parindent}{0.0cm}\textbf{Proof.} By Theorem \ref{TH1}, it is
enough to show that all the edges of $M$ (the unique perfect matching of $%
G_{0}$) are also $\alpha $-critical.\setlength
{\parindent}{3.45ex}

According to Proposition \ref{prop1}, we may write $G$ as $G=S*H$, where $%
S\in \Omega (G)$ and $H=G[V-S]$ has $n(H)=\mu (G)$. By virtue of Lemma \ref
{prop6}($\mathit{ii}$), $G_{0}$ has $\alpha (G_{0})=\left| S-core(G)\right|
=q$ and $core(G_{0})=\emptyset $. Let $M=\{a_{i}b_{i}:1\leq i\leq q\}$ and
suppose that $\{a_{i}:1\leq i\leq q\}=A\subseteq S$. We shall show that any $%
a_{i}b_{i}\in M$ is $\alpha $-critical, by exhibiting a maximum stable set $%
S_{0}$ in $G_{0}$ that satisfies: $b_{i}\in S_{0}$ and $S_{0}\cap
N(a_{i})=\{b_{i}\}$. For the sake of simplicity, let us take $i=1$. In the
sequel, if $D\subseteq V\left( G_{0}\right) $, then by $M(D)$ we mean the
set of vertices, which $D$ is matched onto.

\emph{Claim 1.} There exists some $S_{0}\in \Omega (G_{0})$ with $b_{1}\in
S_{0}$.

Otherwise, any $W\in \Omega (G_{0})$ contains $a_{1}$, because $\left|
M\right| =\left| W\right| $ and $\left| W\cap \{a_{j},b_{j}\}\right| =1$
holds for every $j\in \{1,2,...,q\}$. Hence, it follows that $a_{1}\in
core(G_{0})$, in contradiction with $core(G_{0})=\emptyset $.

\emph{Claim 2.} The following procedure gives rise to some $S_{0}\in \Omega
(G_{0})$ that contains $b_{1}$.\newline

\textbf{Input:} $G_{0},A=\{a_{1},a_{2},...,a_{q}\},b_{1}\in
B=\{b_{1},b_{2},...,b_{q}\}=M(A)$;

\textbf{Output:} $b_{1}\in S_{0}\in \Omega (G_{0})$;

$S_{0}:=\{b_{1}\}$;

$D:=\{b_{1}\}$;

\textbf{while }$(N(D)\cap A)-M(S_{0})\neq \emptyset $ \textbf{do}

\textbf{begin}

\qquad Step 1. $S_{1}:=S_{0}$;

\qquad Step 2. $S_{0}:=S_{0}\cup M((N(D)\cap A)-M(S_{0}))$;

\qquad Step 3. $D:=S_{0}-S_{1}$;

\textbf{end}

Step 4. $S_{0}:=S_{0}\cup M(B-S_{0})$.\newline

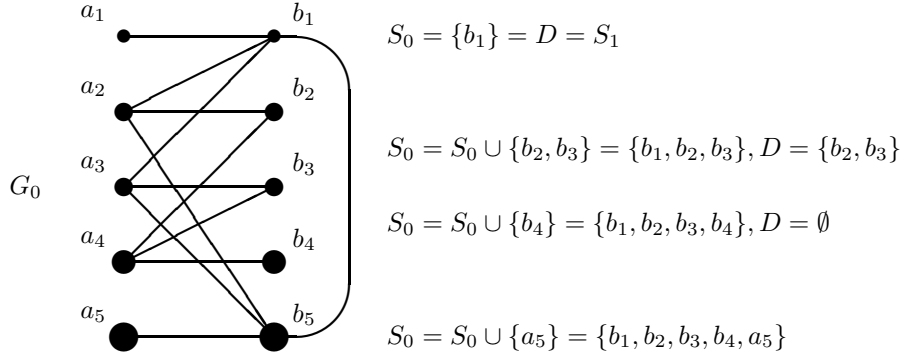
\begin{figure}[h]
\setlength{\unitlength}{1cm}%
\begin{picture}(9,4.5)\thicklines

\put(5.5,0){\makebox(0,0)[l]{$S_{0}=S_{0}\cup \{a_{5}\}=\{b_{1},b_{2},b_{3},b_{4},a_{5}\}$}}
\put(5.5,1.5){\makebox(0,0)[l]{$S_{0}=S_{0}\cup \{b_{4}\}=\{b_{1},b_{2},b_{3},b_{4}\},D=\emptyset $}}
\put(5.5,2.5){\makebox(0,0)[l]{$S_{0}=S_{0}\cup \{b_{2},b_{3}\}=\{b_{1},b_{2},b_{3}\},D=\{b_{2},b_{3}\}$}}
\put(5.5,4){\makebox(0,0)[l]{$S_{0}=\{b_{1}\}=D=S_{1}$}}

  \multiput(2,0)(2,0){2}{\circle*{0.38}}
  \multiput(2,1)(2,0){2}{\circle*{0.32}}
  \multiput(2,2)(2,0){2}{\circle*{0.25}}
  \multiput(2,3)(2,0){2}{\circle*{0.25}}
  \multiput(2,4)(2,0){2}{\circle*{0.19}}

  \put(2,0){\line(1,0){2}}
  \put(2,1){\line(1,0){2}}
  \put(2,2){\line(1,0){2}}
  \put(2,3){\line(1,0){2}}
  \put(2,4){\line(1,0){2}}

  \put(2,2){\line(1,-1){2}} 
  \put(2,1){\line(2,1){2}}
  \put(2,1){\line(1,1){2}}
  \put(2,2){\line(1,1){2}}
  \put(2,3){\line(2,-3){2}}
  \put(2,3){\line(2,1){2}}

   \put(4,2){\oval(2,4)[r]}

\put(1.6,0.3){\makebox(0,0){$a_{5}$}}
\put(4.4,0.3){\makebox(0,0){$b_{5}$}}
\put(1.6,1.3){\makebox(0,0){$a_{4}$}}
\put(4.4,1.3){\makebox(0,0){$b_{4}$}}
\put(1.6,2.3){\makebox(0,0){$a_{3}$}}
\put(4.4,2.3){\makebox(0,0){$b_{3}$}}
\put(1.6,3.3){\makebox(0,0){$a_{2}$}}
\put(4.4,3.3){\makebox(0,0){$b_{2}$}}
\put(1.6,4.3){\makebox(0,0){$a_{1}$}}
\put(4.4,4.3){\makebox(0,0){$b_{1}$}}
\put(0.7,2){\makebox(0,0){$G_{0}$}}

 \end{picture}
\caption{The graph $G_{0}$ has a unique perfect matching and $\xi \left(
G_{0}\right) =0$.}
\label{fig11}
\end{figure}

Clearly, $\left| S_{0}\right| =q$ and no edge of $G_{0}$ joins some $%
a_{l}\in S_{0}$ to any $b_{j}\in S_{0}$, according to building procedure of $%
S_{0}$. Any maximum stable set $W\in \Omega (G_{0})$ that contains $b_{1}$
must contain also all $b_{j}\in S_{0}$, because $\left| W\right| =\left|
M\right| $ and $\left| W\cap \{a_{j},b_{j}\}\right| =1$ holds for every $%
j\in \{1,2,...,q\}$. Hence, the set $\{b_{j}:b_{j}\in S_{0}\}$ is stable,
and consequently, we obtain that $S_{0}\in \Omega (G_{0})$. An example of $%
S_{0}\in \Omega (G_{0})$ obtained by this procedure is illustrated in Figure 
\ref{fig11}. 

\emph{Claim 3.} $S_{0}\cup \{a_{1}\}\in \Omega (G_{0}-a_{1}b_{1})$, and
hence, the edge $a_{1}b_{1}$ is $\alpha $-critical in $G_{0}$.

Firstly, no $a_{i}\in S_{0}$ is adjacent to $a_{1}$, because $a_{i},a_{1}\in
A$. Secondly, no $b_{j}\in S_{0}-\{b_{1}\}$ is adjacent to $a_{1}$,
otherwise there exists an even cycle $C$, with half of its edges belonging
to $M$, which means that $\left( M-E\left( C\right) \right) \cup \left(
E\left( C\right) -M\right) $ is another perfect matching in $G_{0}$, in
contradiction with the premises on $G_{0}$. Therefore, $S_{0}\cup
\{a_{1}\}\in \Omega (G_{0}-a_{1}b_{1})$ and this implies that the edge $%
a_{1}b_{1}$ is $\alpha $-critical in $G_{0}$. Since $a_{1}b_{1}$ is an
arbitrary edge of $M$, we may conclude that all the edges of $M$ are $\alpha 
$-critical in $G_{0}$. \rule{2mm}{2mm}\newline

It is interesting to notice that if $G_{0}$ were bipartite for every
K\"{o}nig-Egerv\'{a}ry graph $G$, then it would be possible to prove Lemma 
\ref{lem3} using only Proposition \ref{prop3}. Figures \ref{fig4}, \ref
{fig11} show that Proposition \ref{prop3} is not enough for our purposes,
because there exist non-bipartite K\"{o}nig-Egerv\'{a}ry graphs $G$ with
nonempty cores and whose $G_{0}=G-N[core(G)]$ have a unique perfect matching.

\begin{theorem}
\label{Th2}Let $G$ be a K\"{o}nig-Egerv\'{a}ry graph and $G_{0}=G-N[core(G)]$%
. Then the following assertions are equivalent:

($\mathit{i}$) $G_{0}$ has a unique perfect matching;

($\mathit{ii}$) $\alpha $-critical edges of $G_{0}$ form a maximal matching
in $G_{0}$;

($\mathit{iii}$) $\xi (G)+\eta (G)=\alpha (G)$;

($\mathit{iv}$) $\sigma (G)+\eta (G)=\mu (G)$;

($\mathit{v}$) $\xi (G)+2\eta (G)+\sigma (G)=n(G)$.
\end{theorem}

\setlength {\parindent}{0.0cm}\textbf{Proof.} According to Proposition \ref
{prop9}, $G_{0}$ is also a K\"{o}nig-Egerv\'{a}ry graph and has a perfect
matching, say $M_{0}$.\setlength
{\parindent}{3.45ex}

($\mathit{i}$) $\Leftrightarrow $ ($\mathit{ii}$) If $M_{0}$ is the unique
perfect matching of $G_{0}$, all its edges are $\mu $-critical and, by Lemma 
\ref{lem3}, $\alpha $-critical, as well. In other words, the $\alpha $%
-critical edges of $G_{0}$ form a maximal matching. The converse is true
according to Proposition \ref{prop5}($\mathit{iii}$).

($\mathit{i}$) $\Rightarrow $ ($\mathit{iii}$) Assume that $M_{0}$ is the
unique perfect matching of $G_{0}$. By Lemma \ref{prop6}($\mathit{ii}$), it
follows that $\alpha (G_{0})=\alpha (G)-\xi (G)$. Lemma \ref{prop6}($\mathit{%
iii}$) and the uniqueness of $M$ imply that $\alpha (G_{0})=\eta
(G_{0})=\eta (G)$. Hence, it results in $\xi (G)+\eta (G)=\alpha (G)$.

($\mathit{iii}$) $\Leftrightarrow $ ($\mathit{iv}$) $\Leftrightarrow $($%
\mathit{v}$) It is the claim of Proposition \ref{prop10}.

($\mathit{v}$) $\Rightarrow $ ($\mathit{ii}$) By Proposition \ref{prop11}, 
\[
|N(core(G))|=|\cap \{V-S:S\in \Omega (G)\}|=\sigma (G). 
\]
Hence, $n(G_{0})=n(G)-\xi (G)-\sigma (G)$. Now, our premise claims that $%
2\eta (G)=n(G_{0})$. By Lemma 3.4($\mathit{iii}$) we obtain $2\eta
(G_{0})=n(G_{0})$. According to Theorem \ref{TH1}($\mathit{iii}$) the set of 
$\alpha $-critical edges of $G$ form a matching, say $M$. Applying again
Lemma 3.4($\mathit{iii}$), we see that $M_{0}=M$ and it consists of $\alpha $%
-critical edges of $G_{0}$. \rule{2mm}{2mm}\newline

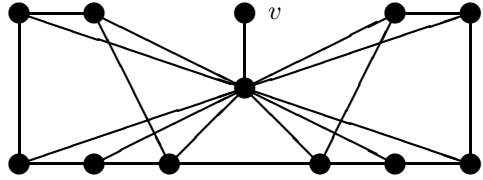
\begin{figure}[h]
\setlength{\unitlength}{1cm}%
\begin{picture}(5,2.5)\thicklines

  \multiput(4,0)(1,0){3}{\circle*{0.29}}
  \multiput(8,0)(1,0){3}{\circle*{0.29}}
  \multiput(4,2)(1,0){2}{\circle*{0.29}}
  \multiput(9,2)(1,0){2}{\circle*{0.29}}
  \multiput(7,1)(0,1){2}{\circle*{0.29}}

  \put(4,0){\line(1,0){6}}
  \put(4,0){\line(0,1){2}} 
  \put(4,2){\line(1,0){1}} 
  
  \put(4,2){\line(3,-1){3}}
  \put(4,0){\line(3,1){3}}
  \put(5,0){\line(2,1){2}}
  \put(6,0){\line(1,1){1}}
  \put(5,2){\line(1,-2){1}}
  \put(5,2){\line(2,-1){2}}
  \put(9,2){\line(1,0){1}}
  \put(10,0){\line(0,1){2}}
  
\put(8,0){\line(1,2){1}}
\put(7,1){\line(1,-1){1}}
\put(7,1){\line(2,1){2}}
\put(7,1){\line(3,1){3}}
\put(7,1){\line(2,-1){2}}
\put(7,1){\line(3,-1){3}}

  \put(7,1){\line(0,1){1}}

\put(7.4,2){\makebox(0,0){$v$}}

 \end{picture}
\caption{A non-Koenig-Egervary graph satisfying $\xi (G)+\eta (G)<\alpha (G)$%
.}
\label{fig5}
\end{figure}

Notice that Theorem \ref{Th2} fails for non-K\"{o}nig-Egerv\'{a}ry graphs.
In Figure \ref{fig5} is presented a non-K\"{o}nig-Egerv\'{a}ry graph $G$
having $\xi (G)=\left| \{v\}\right| =1,\eta (G)=10$, (all the edges of the
two $C_{5}$ are $\alpha $-critical), $\alpha (G)=5<\mu (G)=6$, but $%
G_{0}=G-N[core(G)]$ owns a unique perfect matching.

Now using Theorem \ref{Th2} we are giving a new characterization of the
bipartite graphs that have a unique perfect matching (see some previous
discussions of this topic in \cite{Cechlarova} and \cite{Martinez}). This
result generalizes Corollary \ref{cor3}.

\begin{corollary}
\label{cor6}Let $G$ be a bipartite graph. Then the following assertions are
equivalent:

($\mathit{i}$) $G$ has a unique perfect matching;

($\mathit{ii}$) $\alpha $-critical edges of $G$ form a maximal matching;

($\mathit{iii}$) $\eta (G)=\alpha (G)$;

($\mathit{iv}$) $\eta (G)=\mu (G)$;

($\mathit{v}$) $2\eta (G)=n(G)$.
\end{corollary}

\setlength {\parindent}{0.0cm}\textbf{Proof.} ($\mathit{i}$) $%
\Leftrightarrow $ ($\mathit{ii}$) If $M$ is the unique perfect matching of $%
G $, all its edges are $\mu $-critical and, by Proposition \ref{prop3}, $%
\alpha $-critical, as well. In other words, the $\alpha $-critical edges of $%
G$ form a maximal matching. The converse is true according to Proposition 
\ref{prop5}($\mathit{iii}$).\setlength
{\parindent}{3.45ex}

The other equivalences follow from Theorem \ref{Th2}, and the observation
that if a bipartite graph has a perfect matching, then the two stable sets
of its standard partition are maximum, and, consequently, $\xi (G)=0$. \rule%
{2mm}{2mm}\newline

It is interesting to notice that the equality $2\alpha (G)=n(G)$ mentioned
in Corollary \ref{cor3} follows from Corollary \ref{cor6}, but it can not
join the above series of equivalences (see, for example, $C_{4}$).

\begin{figure}[h]
\setlength{\unitlength}{1.0cm}%
\begin{picture}(5,1.1)\thicklines
  \multiput(5,0)(1,0){4}{\circle*{0.29}}
  \multiput(6,1)(1,0){3}{\circle*{0.29}}
  \put(5,0){\line(1,0){3}}
  \put(7,1){\line(1,0){1}}
  \multiput(6,0)(1,0){3}{\line(0,1){1}}
 \end{picture}
\caption{$\xi (G)=2,\eta (G)=0,\alpha (G)=4,\sigma (G)=1,\mu (G)=4,n(G)=7$.}
\label{fig33}
\end{figure}
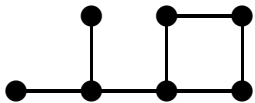

Let us also observe that for the bipartite graph $G$ in Figure \ref{fig33},
the subgraph $G_{0}=G-N[core(G)]$ has more than one perfect matching.

\begin{proposition}
\label{prop4}If $G$ is a K\"{o}nig-Egerv\'{a}ry graph and there is some $%
S\in \Omega (G)$ such that the set $W=(S,V(G)-S)$ generates a forest, then 
\[
\xi (G)+\eta (G)=\alpha (G),\sigma (G)+\eta (G)=\mu (G),\ and\ \xi (G)+2\eta
(G)+\sigma (G)=n(G).
\]
\end{proposition}

\setlength {\parindent}{0.0cm}\textbf{Proof.} If $%
G_{0}=G-N[core(G)],A=S-core(G),B=V(G)-S-N(core(G))$, then Proposition \ref
{prop9}($\mathit{iii}$) implies that $G_{0}$ is also a
K\"{o}nig-Egerv\'{a}ry graph and has a perfect matching, say $M$. Let $G_{1}$
be the partial graph of $G_{0}$ having $W\cap E(G_{0})$ as edge set. Then, $%
M $ is a perfect matching in $G_{1}$, as well. Since $G_{1}$ is a forest, $M$
is unique. By Lemma \ref{lem2}, any maximum matching of $G_{0}$ is contained
in $(A,B)$, and since the edges from $(A,B)$ yield a unique perfect
matching, namely $M$, it follows that $M$ is the unique perfect matching of $%
G_{0}$ itself. Hence, according to Theorem \ref{Th2}, we obtain that $\xi
(G)+\eta (G)=\alpha (G)$. By Proposition \ref{prop7}($\mathit{iii}$), it
implies $\sigma (G)+\eta (G)=\mu (G)$, and immediately $\xi (G)+2\eta
(G)+\sigma (G)=n(G)$. \rule{2mm}{2mm}\setlength
{\parindent}{3.45ex}\newline

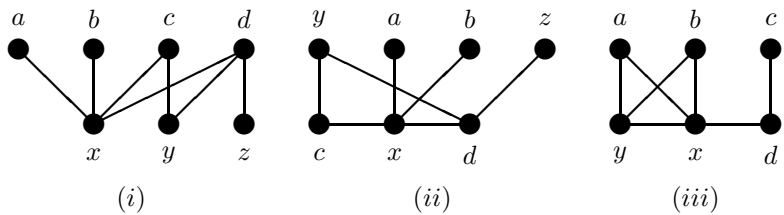
\begin{figure}[h]
\setlength{\unitlength}{1cm}%
\begin{picture}(5,2.5)\thicklines

  \multiput(2,1)(1,0){3}{\circle*{0.29}}
  \multiput(1,2)(1,0){4}{\circle*{0.29}}
  \put(2,1){\line(0,1){1}}
  \put(3,1){\line(0,1){1}} 
  \put(4,1){\line(0,1){1}} 
  
  \put(1,2){\line(1,-1){1}}
  \put(2,1){\line(1,1){1}}
  \put(2,1){\line(2,1){2}}
  \put(3,1){\line(1,1){1}}
 
\put(1,2.4){\makebox(0,0){$a$}}
\put(2,2.4){\makebox(0,0){$b$}}
\put(3,2.4){\makebox(0,0){$c$}}
\put(4,2.4){\makebox(0,0){$d$}}
\put(2,0.6){\makebox(0,0){$x$}}
\put(3,0.6){\makebox(0,0){$y$}}
\put(4,0.6){\makebox(0,0){$z$}}
\put(2.5,0){\makebox(0,0){$(i)$}}  

  \multiput(5,2)(1,0){4}{\circle*{0.29}}
  \multiput(5,1)(1,0){3}{\circle*{0.29}}
  \put(5,1){\line(0,1){1}}
  \put(6,1){\line(0,1){1}}
  \put(5,1){\line(1,0){2}}
  \put(5,2){\line(2,-1){2}}
  \put(6,1){\line(1,1){1}}
  \put(7,1){\line(1,1){1}}
  
\put(5,2.4){\makebox(0,0){$y$}}
\put(6,2.4){\makebox(0,0){$a$}}
\put(7,2.4){\makebox(0,0){$b$}}
\put(8,2.4){\makebox(0,0){$z$}}
\put(5,0.6){\makebox(0,0){$c$}}
\put(6,0.6){\makebox(0,0){$x$}}
\put(7,0.6){\makebox(0,0){$d$}}
\put(6.5,0){\makebox(0,0){$(ii)$}}  

  \multiput(9,1)(1,0){3}{\circle*{0.29}}
  \multiput(9,2)(1,0){3}{\circle*{0.29}}
  \put(9,1){\line(0,1){1}}
  \put(10,1){\line(0,1){1}} 
  \put(11,1){\line(0,1){1}} 
  \put(9,1){\line(1,0){2}}
  \put(9,1){\line(1,1){1}}
  \put(9,2){\line(1,-1){1}}
 
\put(9,2.4){\makebox(0,0){$a$}}
\put(10,2.4){\makebox(0,0){$b$}}
\put(11,2.4){\makebox(0,0){$c$}}
\put(11,0.6){\makebox(0,0){$d$}}
\put(10,0.6){\makebox(0,0){$x$}}
\put(9,0.6){\makebox(0,0){$y$}}
\put(10,0){\makebox(0,0){$(iii)$}}  

 \end{picture}
\caption{K\"{o}nig-Egervary graphs satisfying $\xi (G)+\eta (G)=\alpha (G)$.}
\label{fig6}
\end{figure}

It is worth observing that if $(S,V(G)-S)$ generates a forest for some $S\in
\Omega (G)$, this is not necessarily true for all maximum stable sets of $G$%
. For example, the graph $G$ presented in Figure \ref{fig6}($\mathit{i}$)
and Figure \ref{fig6}($\mathit{ii}$) has the partition $\{S_{1}=\{a,b,c,d\}%
\in \Omega (G),V(G)-S_{1}\}$ such that $(S_{1},V(G)-S_{1})$ does not
generate a forest, (see Figure \ref{fig6}($\mathit{i}$)), while for the
partition $\{S_{2}=\{a,b,y,z\}\in \Omega (G),V(G)-S_{2}\}$ the set $%
(S_{2},V(G)-S_{2})$ generates a forest (see Figure \ref{fig6}($\mathit{ii}$%
)). Let us also remark that the converse of Proposition \ref{prop4} is not
generally true. For instance, the graph in Figure \ref{fig6}($\mathit{iii}$)
is a counterexample. 

\begin{corollary}
\label{cor1}If $T$ is a tree, then 
\[
\xi (T)+\eta (T)=\alpha (T),\sigma (T)+\eta (T)=\mu (T),\ and\ \xi (T)+2\eta
(T)+\sigma (T)=n(T).
\]
\end{corollary}

As a consequence of Corollary \ref{cor1}, we obtain:

\begin{corollary}
\cite{Zito}\label{cor5} If $T$ is a tree, then a vertex $v\in V(T)$ is in
some but not in all maximum stable sets of $T$ if and only if $v$ is an
endpoint of an $\alpha $-critical edge.
\end{corollary}

\setlength {\parindent}{0.0cm}\textbf{Proof.} If $v\in V(T)$ is in some but
not in all maximum stable sets of $T$, then there exists $S\in \Omega (T)$
such that $v\in S-core(T)$. By Theorem \ref{TH1}, $\alpha $-critical edges
of $T$ form a matching. Proposition \ref{prop3} ensures that they are also $%
\mu $-critical, because $T$ is bipartite. Consequently, these edges belong
to any maximum matching, which, according to Lemma \ref{lem2}, is included
in $(S,V(T)-S)$. Since, by Lemma \ref{prop6}($\mathit{i}$), no $\alpha $%
-critical edge has an endpoint in $N[core(T)]$, and Corollary \ref{cor1}
ensures that $\eta (T)=\alpha (T)-\xi (T)=|S-core(T)|$, we infer that $v$
must be an endpoint of an $\alpha $-critical edge.%
\setlength
{\parindent}{3.45ex}

Conversely, let $e=vw$ be an $\alpha $-critical edge in $T$ and $\overline{S}%
\in \Omega (T-e)$. Since $\left| \overline{S}\right| =\alpha (T-e)>\alpha
(T) $, it follows that $v,w\in \overline{S}$ and therefore, $\overline{S}%
-\{v\},\overline{S}-\{w\}\in \Omega (T)$. Hence, $v$ is in some, namely, in $%
\overline{S}-\{w\}$, but not in all maximum stable sets of $T$, namely, not
in $\overline{S}-\{v\}$. \rule{2mm}{2mm}\newline

Notice that Corollary \ref{cor1} and Corollary \ref{cor5} are not valid for
general bipartite graphs (see, for instance, the graph in Figure \ref{fig33}%
).

\section{Conclusions}

In this paper we state several properties of $\alpha $-critical and $\mu $%
-critical edges belonging to K\"{o}nig-Egerv\'{a}ry graphs. These findings
generalize some previously known results for trees and bipartite graphs. We
have proved that for bipartite graphs and for some special
K\"{o}nig-Egerv\'{a}ry graphs, their sets of $\alpha $-critical edges and $%
\mu $-critical edges coincide. It seems to be interesting to characterize
all the graphs having this property. From the other point of view, since the 
$\alpha $-critical edges of a K\"{o}nig-Egerv\'{a}ry graph span disjoint
cliques of order two, one may be interested in describing the type of graphs
where their $\alpha $-critical edges span disjoint cliques of order larger
than two. Another challenging problem is to describe classes of
non-K\"{o}nig-Egerv\'{a}ry graphs $G$ satisfying $\xi (G)+\eta (G)=\alpha
(G) $, $\xi (G)+\eta (G)\leq \alpha (G)$, and/or $\alpha (G)+\sigma (G)=\mu
(G)+\xi (G)$.

\end{document}